\documentclass[a4paper,10pt]{article}
\usepackage{bm,amsmath,amssymb,mathrsfs,graphicx,citesort,enumerate}
\newcommand{\nc}{\newcommand}
\nc{\rnc}{\renewcommand}
\nc{\beq}{\begin{equation}}
\nc{\eeq}{\end{equation}}
\nc{\nn}{\nonumber}
\rnc{\(}{\left(}
\rnc{\)}{\right)}
\rnc{\[}{\left[}
\rnc{\]}{\right]}

\newtheorem{definition}{Definition}
\newtheorem{theorem}{Theorem}
\newtheorem{proposition}{Proposition}
\newtheorem{lemma}{Lemma}

\numberwithin{definition}{section}
\numberwithin{equation}{section}
\numberwithin{lemma}{section}
\numberwithin{proposition}{section}
\numberwithin{theorem}{section}
\numberwithin{corollary}{section}

\nc{\sq}{\qquad $\blacksquare$}
\nc{\wsq}{\qquad $\square$}

\begin{document}%
%
\title{Nonstandard representations of type $C$ affine Hecke algebra from
$K$-operators}
\author{Kohei Motegi\thanks{E-mail: motegi@gokutan.c.u-tokyo.ac.jp} \\
\it Okayama Institute for Quantum Physics, \\
\it Kyoyama 1-9-1,
Okayama 700-0015, Japan, 
\\}

\date{\today}
 

\maketitle

%
%
\begin{abstract}
We construct nonstandard finite-dimensional
representations of type $C$ affine Hecke algebra
from the viewpoint of quantum integrable models.
There exists two classes of nonstandard solutions to the Yang-Baxter equation
called the Cremmer-Gervais and Jordanian $R$-matrices.
These $R$-matrices also satisfy the Hecke-relation, thus can be used to
construct nonstandard finite-dimensional representations of
type $A$ affine Hecke algebra.
We construct the corresponding nonstandard representations for
type $C$ affine Hecke algebra by explicitly constructing solutions to the
reflection equation under the Hecke relation.
We achieve this by taking the finite-dimensional representations
and deBaxterizing the $K$-operators acting on the
infinite-dimensional function space, taking advantage of the fact that
the Cremmer-Gervais and Jordanian $R$-matrices can be obtained from
the $R$-operator.
\end{abstract}

Mathematics Subject Classification (2010). 20C08, 16T25, 14H70, 17B37. \\
\\

Keywords. Affine Hecke algebra, Yang-Baxter equation, reflection equation,
quantum integrable models.

%
\section{Introduction}
Affine Hecke algebras plays special roles in
mathematics and mathematical physics.
They are not only one of the most important algebras in
representation theory related to Yangians and quantum affine algebras
for example, but also have various
applications to other branches from mathematical physics,
knot theory to the recent progress in categorification,
geometric representation theory and so on (see \cite{Ma,Ch} for example for
general treatments and applications of affine Hecke algebras).

From the point of view of mathematical physics,
affine Hecke algebras are very intimate with quantum integrable models,
in particular with the object called the $R$-matrix.
Quantum integrable models are special classes of quantum
many-body systems equipped with special algebraic structures.
At the time when the very first example of
quantum integrable model was constructed \cite{Bethe},
Bethe treated only the Hamiltonian,
but it is in general very difficult to judge whether
a Hamiltonian is quantum integrable or not.
The advances in the study of quantum integrable models
lead to the discovery of the Yang-Baxter relation
which is a characteristic feature for quantum integrable models,
and it is widely believed today that
the $R$-matrix satisfying the Yang-Baxter relation
is the most fundamental object in quantum integrable models
since the global commuting transfer matrix, which is a generating function
of conserved quantities including Hamiltonian, is constructed
from the local $R$-matrix.

Based on the $R$-matrix, there are various algebraic methods
developed to analyze quantum integrable models and algebraic structures
investigated from the point of view of quantum integrability.
(see the papers and books
\cite{FST,Ba,KBI,JM} about quantum inverse scattering method
and $q$-vertex operator).
The relation between $R$-matrices and affine Hecke algebras
is one such example.
It is well-known that 
trigonometric $R$-matrices can be constructed from affine Hecke algebra
of type $A$, and vice versa. In fact, the standard $R$-matrix of the
quantum group $U_q(\widehat{sl_2})$ \cite{Dr,Ji} without spectral parameter
is nothing but a representation for a generator of the affine Hecke algebra.
The braid relation which is one type of the defining relations 
among the generators for the affine Hecke algebra
is nothing but the Yang-Baxter relation in quantum integrable models.
The trigonometric $R$-matrix with spectral parameter
can be constructed as a combination of a generator and a permutation operator.
This construction can be generalized to other root systems.

For the affine Hecke algebra of type $C$,
the boundary condition of the corresponding quantum integrable models
is modified from the standard periodic boundary condition to
the open boundary condition.
From the point of view of quantum integrable models,
two of the generators of type $C$ affine Hecke algebra
are nothing but the solutions to the reflection equation \cite{Sklyanin,KS},
which is the boundary condition on quantum systems
to ensure integrability of the model.
There are many solutions to the reflection equation
constructed today (see \cite{DG,IK,Ya,KM} for example).

In this paper, we construct nonstandard finite-dimensional representations
of type $C$ affine Hecke algebra.
Despite its extensive studies on various aspects,
it seems that the problem of concretely realizing representations of
affine Hecke algebra is not well investigated.
We approach this problem by using the power of quantum integrable models.
From the point of view of affine Hecke algebra,
the interesting $R$-matrices are those satisfying
the Hecke relation.
Besides the standard $R$-matrix of the quantum group $U_q(\widehat{sl_2})$,
there is a nonstandard one called the Cremmer-Gervais (CG) $R$-matrix.
This $R$-matrix originally appeared in the context of the
Toda field theory \cite{CG} as a constant $R$-matrix,
and the Baxterized $R$-matrix was derived in \cite{EH}.
The CG-$R$-matrix satisfies the Hecke relation,
thus can be served as a representation for type $A$ affine Hecke algebra.
We call the corresponding
representation as a nonstandard trigonometric representation
since it comes from the nonstandard $R$-matrix.

In this paper, we construct nonstandard
representations for type $C$ affine Hecke algebra.
From the point of view of quantum integrable models,
the problem of constructing nonstandard trigonometric representations
for type $C$ is equivalent to finding solutions to the reflection equation
of the CG-$R$-matrix under the Hecke relation.

To achieve this, we first reivew how the Baxterized
CG-$R$-matrix was derived \cite{EH}.
It was derived by taking finite-dimensional representation of the
trigonometric limit of the elliptic Shibukawa-Ueno (SU) $R$-operator \cite{SU},
which is an infinite-dimensional $R$-operator acting on the space of functions.
We remark here
that taking the finite-dimensional representation of the elliptic
SU-$R$-operator gives the Baxter-Belavin model
\cite{Baxtereightvertex,Belavin},
whose trigonometric limit is the standard $R$-matrix of $U_q(\widehat{sl_2})$.
The nonstandard CG-$R$-matrix is obtained in a different manner.
We first take the trigonometric limit of the SU-$R$-operator.
Taking the finite-dimensional representation
in the trigonometric basis after that procedure
gives the CG-$R$-matrix
(see also \cite{AHZ} for a similar construction).

By taking advantage of this fact, we construct solutions to the
reflection equation of the CG-$R$-matrix in the following way.
We start from the elliptic $K$-operator by Hikami-Komori (HK) \cite{Hi,HK}
which are solutions to the reflection equation of
the SU-$R$-operator.
We first take the trigonometric limit of the elliptic $K$-operator
to get the trigonometric $K$-operator.
We next take the finite-dimensional representation in the trigonometric
basis to get the Baxterized $K$-matrix.
By construction, the $K$-matrix is the solution to the
reflection equation of the Baxterized $R$-matrix.
To extract representations of the affine Hecke algebra,
we need one more thing to do. Namely, one has to do the
deBaxterization, i.e., extract the constant
deBaxterized $K$-matrix out of the spectral parameter dependent $K$-matrix
by getting rid of the spectral parameter.
We show from its construction that the $K$-matrix
satisfy both the constant reflection equation and the Hecke relation,
which are exactly the defining relations the boundary generators
of affine Hecke algebra of type $C$ must satisfy.

We also construct another type of nonstandard representation for
a special case of affine Hecke algebra of type $C$
by considering a class of quantum integrable models of rational type.
The rational Jordanian $R$-matrix and its assoiciated $K$-matrix
serve as representations for the generators.
We also find the representations in a similar but different way.
We first go furthermore to the rational limit
of the trigonometric $R$-operator and the $K$-operator,
next take finite-dimensional representations
in the rational basis.
Then we deBaxterize to get the $R$-matrix and its corresponding
$K$-matrix.

In the next section, we state the main results about
nonstandard representations of type $C$ affine Hecke algebra,
and give the outline of the proof.
In section 3, we first review the CG-$R$-matrix
and its construction from the SU-$R$-operator
by taking the finite-dimensional representation in the trigonometric basis
and deBaxterizing.
Then we apply the same degeneration procedure to construct
solutions to the reflection equation to the CG-$R$-matrix
from the HK-$K$-operator,
and deBaxterize to get representations
for the boundary generators of type $C$ affine Hecke algebra,
which gives the proof for one of the main results.
In section 4, we apply a similar procedure to the rational case to
give the other main results.
Section 5 is devoted to the conclusion.

\section{Type $C$ affine Hecke algebra and nonstandard representations}
In this paper, we construct representations of the
following type $C$ affine Hecke algebra.
\begin{definition}
Type $C$ affine Hecke algebra $H_n=H_n(t,t_n,t_0)$
is defined as an algebra generated by $T_j, \ j=0,\cdots,n$
satisfying the following relations
\begin{align}
&(T_0-t_0)(T_0+t_0^{-1})=0, \label{relone} \\
&(T_j-t)(T_j+t^{-1})=0, \ \ \ 1 \le j \le n-1, \label{reltwo} \\
&(T_n-t_n)(T_n+t_n^{-1})=0, \label{relthree} \\
&T_0 T_1 T_0 T_1=T_1 T_0 T_1 T_0, \label{relfour} \\
&T_j T_{j+1} T_j=T_{j+1} T_j T_{j+1}, \ \ \ 1 \le j \le n-2, \label{relfive} \\
&T_{n-1}T_n T_{n-1}T_n=T_n T_{n-1}T_n T_{n-1}, \label{relsix} \\
&T_j T_k=T_k T_j, \ \ \ |j-k| \ge 2. \label{relseven}
\end{align}
\hspace*{\fill}    \sq
\end{definition}
We give nonstandard representations
for type $C$ affine Hecke algebra by extracting the generators
from a class of trigonometric and rational quantum integrable models.
To state the main theorem, we first fix notations.
Let $V$ be an $N$-dimensional complex vector space.
We denote the orthonormal basis of $V$ by $\{ e_j, j=0,\cdots,N-1 \}$.
The matrix element $[A]_{j}^k$ of $A \in \mathrm{End}(V)$
with respect to this basis is defined as
\begin{eqnarray}
A e_j=\sum_{k=0}^{N-1} e_k [A]_{j}^k. \nonumber \\
\end{eqnarray}
The permutation matrix $P$ is defined as
\begin{eqnarray}
P(x \otimes y)=y \otimes x \ \mathrm{for \  any} \ x,y \in V.
\end{eqnarray}
Let $V^{\otimes n}=V_1 \otimes V_2 \otimes \cdots \otimes V_n$
be a tensor product of complex vector spaces.
For a matrix $A \in V$, let us define
a matrix $A_j \in \mathrm{End} (V_1 \otimes \cdots \otimes V_n)$ as
a matrix acting on the complex vector space $V_j$ as $A$,
and acting on the remaining
complex vector spaces $V_k \ (k \neq j)$ as an identity matrix.
For a matrix $A_{jk} \in \mathrm{End} (V_j \otimes V_k)$,
we define $\check{A}_{jk}$ as
\begin{eqnarray}
\check{A}_{jk}=A_{jk}P_{jk}.
\end{eqnarray}
\begin{definition} \label{definitionrkmatrixone}
We define the $R$-matrix $R^{\mathrm{tr}}(q,p) \in V \otimes V$
and the $K$-matrix $K^{\mathrm{tr}}(r,s) \in V$ as
\begin{eqnarray}
\[ R^{\mathrm{tr}}(q,p) \]_{ij}^{kl}
=p^{2(j-k)} \times
 \left\{
\begin{array}{cc}
q, & \mathrm{for} \ i=k \ge j=l, \\
q^{-1}, & \mathrm{for} \ i=k <  j=l, \\
-q+q^{-1}, & \mathrm{for} \ i<k<j, \ i+j=k+l, \\
q-q^{-1}, & \mathrm{for} \ j \le k<i, \
i+j=k+l, \\
0, & \mathrm{otherwise}. \label{cremmer-gervais-r}
\end{array}
\right.
\end{eqnarray}
\begin{eqnarray}
\[ K^{\mathrm{tr}}(r,s) \]_{j}^{k}
=s^{j-k} \times
 \left\{
\begin{array}{cc}
-r^{-1}, & \mathrm{for} \ j \le k, \ j+k=N-1, \\
-r, & \mathrm{for} \ j>k, \ j+k=N-1, \\
r-r^{-1}, & \mathrm{for} \ j \le k<N-1-j, \\
-r+r^{-1}, & \mathrm{for} \ N-1-j<k<j, \\
0, & \mathrm{otherwise}. \label{cremmer-gervais-k}
\end{array}
\right.
\end{eqnarray}
\hspace*{\fill}    \sq
\end{definition}
Using the matrices defined above,
we have the following representation for
type $C$ affine Hecke algebra.
\begin{theorem} (Nonstandard trigonometric representation) \label{thone}
Let $\hat{T}_j, \ j=0 \cdots,n$ be the following matrices
acting on the spin module $V^{\otimes n}$
\begin{align}
&\hat{T}_0=K^{\mathrm{tr}}_1(t_0,s_0), \\
&\hat{T}_j=\check{R}^{\mathrm{tr}}_{j,j+1}(t,p), \ j=1,\cdots,n-1, \\
&\hat{T}_n=K^{\mathrm{tr}}_n(t_n,s_n). \\
\end{align}
The map $\rho$ defined as
\begin{align}
\rho(T_j)=\hat{T}_j, \ j=0,\cdots,n,
\end{align}
is a representation map $\rho:H_n(t,t_n,t_0)
\longrightarrow \mathrm{End}(V^{\otimes n})$.
Namely, 
$\{ \hat{T}_j, \ j=0 \cdots,n \}$
gives a representation for type $C$ affine Hecke algebra
$H_n(t,t_n,t_0)$.
\hspace*{\fill}    \sq
\end{theorem}
We call this representation as a nonstandard trigonometric representation
of type $C$ affine Hecke algebra.
The term ``nonstandard trigonometric" means that this representation
differs from the finite-dimensional representations
based on the standard trigonometric
$R$-matrix of the quantum group $U_q(\widehat{sl_2})$ by Drinfeld and Jimbo,
and its associated $K$-matrix.
``Nonstandard trigonometric" means that this representation
comes from the nonstandard trigonometric $R$-matrix
called the Cremmer-Gervais (CG) $R$-matrix and its associated $K$-matrix
defined in \eqref{cremmer-gervais-r} and \eqref{cremmer-gervais-k}
which we show in this paper.

For the case when the parameters of the type $C$ affine Hecke algebra
is special $t=t_0=t_n=1$,
we can also construct another nonstandard representation.
\begin{definition} \label{definitionrkmatrixtwo}
We define the $R$-matrix $R^{\mathrm{ra}}(\kappa,h) \in V \otimes V$
and the $K$-matrix $K^{\mathrm{ra}}(\nu,g) \in V$ as
\begin{align}
&\[ R^{\mathrm{ra}}(\kappa,h) \]_{ij}^{kl}
\nonumber \\
=&(-1)^{j-l} h^{i+j-k-l} \Bigg\{
\left(
\begin{array}{c}
i \\
k
\end{array}
\right)
\left(
\begin{array}{c}
j \\
l
\end{array}
\right)
-\frac{\kappa}{h}
\sum_m
(-1)^{m-k}
\left(
\begin{array}{c}
i \\
m
\end{array}
\right)
\left(
\begin{array}{c}
j+m-k-1 \\
l
\end{array}
\right)
\epsilon(j,m,k)
\Bigg\},
\label{jordanian-r}
\end{align}
where $\epsilon(i,j,k)$ is defined as
\begin{eqnarray}
\epsilon(i,j,k)
=\left\{
\begin{array}{cc}
1, & \mathrm{for} \ i \le k < j, \\
-1, & \mathrm{for} \ j \le k < i, \\
0, & \mathrm{otherwise}.
\end{array}
\right.
\end{eqnarray}
\begin{eqnarray}
\[ K^{\mathrm{ra}}(\nu,g) \]_{j}^{k}
=(-1)^j
\left(
\begin{array}{c}
j \\
k
\end{array}
\right)
g^{j-k}
+2 \nu \sum_{0 \le l<j}(-1)^{j-l}
\left(
\begin{array}{c}
j-l-1 \\
k-l
\end{array}
\right)
g^{j-k-1}
.
\label{jordanian-k}
\end{eqnarray}
\hspace*{\fill}    \sq
\end{definition}
The following serves as another representation for
type $C$ affine Hecke algebra when $t=t_0=t_n=1$.
\begin{theorem} (Nonstandard rational representation) \label{thtwo}
Let $\hat{T}_j, \ j=0 \cdots,n$ be the following matrices
acting on the spin module $V^{\otimes n}$
\begin{align}
&\hat{T}_0=K^{\mathrm{ra}}_1(\nu_0,g_0), \\
&\hat{T}_j=\check{R}^{\mathrm{ra}}_{j,j+1}(\kappa,h), \ j=1,\cdots,n-1, \\
&\hat{T}_n=K^{\mathrm{ra}}_n(\nu_n,g_n). \\
\end{align}
The map $\rho$ defined as
\begin{align}
\rho(T_j)=\hat{T}_j, \ j=0,\cdots,n,
\end{align}
is a representation map $\rho:H_n(1,1,1) 
\longrightarrow \mathrm{End}(V^{\otimes n})$.
Namely, $\{ \hat{T}_j, \ j=0 \cdots,n \}$
gives a representation for type $C$ affine Hecke algebra
$H_n(1,1,1)$.
\hspace*{\fill}    \sq
\end{theorem}
We call the above representation as ``nonstandard rational" representation
for type $C$ affine Hecke algebra
since the generators are coming from a class of rational solutions
\eqref{jordanian-r}
to the Yang-Baxter relation
called the Jordanian $R$-matrix,
and its associated $K$-matrix \eqref{jordanian-k}
which satisfies the corresponding reflection equation.

\section{Outline of the proof}
In this section, we give the proof
of the main theorems.
We explain the procedures of the degeneration and deBaxterization
to extract nonstandard representations
from quantum integrable models.
We give the outline in this section,
and give the detailed lemmas and propositions
essential to construct
trigonometric and rational representations
in sections 4 and 5, respectively.
\subsection{Trigonometric representation}
We construct a representation of type $C$ affine Hecke algebra
by extracting the generators
from a class of trigonometric quantum integrable models.
In quantum integrable models,
the $R$-matrix $R(\lambda)$ satisfying the Yang-Baxter relation
\begin{equation}
R_{12}(\lambda)R_{13}(\lambda+\mu)R_{23}(\mu)=
R_{23}(\mu)R_{13}(\lambda+\mu)R_{12}(\lambda), \label{yang-baxter}
\end{equation}
is the fundamental object.
The parameter $\lambda$ of $R(\lambda)$ is called the spectral parameter,
and is important to treat quantum integrable models.
However, to construct representations of affine Hecke algebras,
we want to get rid of it.
If an $R$-matrix $R(\lambda)$ satisfying the Yang-Baxter relation
\eqref{yang-baxter} can be decomposed
using the permutation $P$ and the $\lambda$-independent $R$ as
\begin{equation}
R(\lambda)=f(\lambda)(P+g(\lambda)R), \ g(\lambda)=
\frac{\mathrm{e}^{-2 \pi \mathrm{i} \lambda}-1}{q-q^{-1}}, \label{baxterone}
\end{equation}
and the $\lambda$-independent $R$ satisfies the Hecke relation
\begin{equation}
(\check{R}-q)(\check{R}+q^{-1})=0, \label{baxtertwo}
\end{equation}
where $\check{R}=RP$,
then one can show that the $\lambda$-independent
$R$ satisfies the Yang-Baxter relation
without spectral parameter
\begin{equation}
R_{12}R_{13}R_{23}=
R_{23}R_{13}R_{12},
\end{equation}
which is equivalent to the braid relation for $\check{R}=RP$
\begin{equation}
\check{R}_{12}\check{R}_{23}\check{R}_{12}
=\check{R}_{23}\check{R}_{12}\check{R}_{23}.
\end{equation}
This shows that $\check{R}_{j,j+1}, \ j=1,\cdots,n-1$ can be used as a
representation for the generators $T_j, \ j=1,\cdots,n-1$.
The braid relation
\begin{equation}
\check{R}_{j,j+1}\check{R}_{j+1,j+2}\check{R}_{j,j+1}
=\check{R}_{j+1,j+2}\check{R}_{j,j+1}\check{R}_{j+1,j+2}, \ 1 \le j \le n-2,
\end{equation}
corresponds to a defining relation \eqref{relfive}
for the type $C$ affine Hecke algebra,
and the Hecke relation \eqref{reltwo}
can be represented by the Hecke relation for the matrix $\check{R}$
under the identification of parameters $q=t$
\begin{equation}
(\check{R}_{j,j+1}-q)(\check{R}_{j,j+1}+q^{-1})=0.
\end{equation}
We refer to this procedure to construct representations
from quantum integrable models as deBaxterization.
The CG-$R$-matrix realizes a deBaxterization procedure
to construct representations for type $A$ affine Hecke algebra,
which we explain later.

If once a representation for $T_j, \ j=1,\cdots,n-1$ is constructed,
the remaining thing is to construct
the representations for the generators $T_0$, $T_n$
which satisfy the defining relations
\eqref{relone}, \eqref{relthree}, \eqref{relfour} and \eqref{relsix}.
The representation can also be constructed from objects of
quantum integrable models.
The quantum integrability under the reflecting boundary condition
is ensured by the reflection equation
\begin{align}
&R_{21}(\lambda_{1}-\lambda_{2})K_{1}(\lambda_{1})
R_{12}(\lambda_{1}+\lambda_{2})K_{2}(\lambda_{2})
=K_{2}(\lambda_{2})R_{21}(\lambda_{1}+\lambda_{2})
K_{1}(\lambda_{1})R_{12}(\lambda_{1}-\lambda_{2}),
\label{reflection}
\end{align}
where $R_{21}(\lambda)=P_{12}R_{12}(\lambda)P_{12}$.
The representations for generators $T_0$ and $T_{n}$
can be constructed from the $K$-matrix $K(\lambda)$
satisfying the reflection equation \eqref{reflection}
in the following way.

If $K(\lambda)$ satisfying the reflection equation \eqref{reflection}
can be decomposed using the identity and the $\lambda$-independent $K$ as
\begin{equation}
K(\lambda)=a(\lambda)(I+b(\lambda)K), \ b(\lambda)
=\frac{\mathrm{e}^{-4 \pi \mathrm{i} \lambda}-1}{r-r^{-1}}, \label{baxterthree}
\end{equation}
and the $\lambda$-independent $K$ satisfies the Hecke relation
\begin{equation}
(K-r)(K+r^{-1})=0, \label{baxterfour}
\end{equation}
then $K$ satisfies the reflection equations
\begin{align}
&\check{R}_{12}K_{1}\check{R}_{12}K_{1}
=K_{1}\check{R}_{12}K_{1}\check{R}_{12}, \\
&\check{R}_{n-1,n}K_{n}\check{R}_{n-1,n}K_{n}
=K_{n}\check{R}_{n-1,n}K_{n}\check{R}_{n-1,n},
\end{align}
which can be identified as representations for
the defining relations
\eqref{relfour} and \eqref{relsix}.
The Hecke relations
\begin{align}
&(K_1-r_0)(K_1+r_0^{-1})=0, \\
&(K_n-r_n)(K_n+r_n^{-1})=0,
\end{align}
can be identified with the Hecke relations
for $T_0$ and $T_n$
\eqref{relone} and \eqref{relthree}
under the identification of boundary parameters $r_0=t_0$ and $r_n=t_n$.
Thus, the $K$-matrices $K_1$ and $K_n$
serves as representations for $T_0$ and $T_n$ respectively.
We call this procedure to construct representations
of $T_0$ and $T_n$ from the reflection equation \eqref{reflection}
as deBaxterization for boundary.

The above is a description of
a general procedure to construct representations of
type $C$ affine Hecke algebra $H_n(t,t_n,t_0)$
from quantum integrable models by deBaxterization to get rid of
the spectral parameters.
This procedure applies
if the extracted $\lambda$-independent $R$-matrix $R$ and
$K$-matrix $K$ satisfy the Hecke relations.
Now we construct a representation from a class of
trigonometric quantum integrable models.
We start from the Cremmer-Gervais $R$-matrix $R^{\mathrm{tr}}(q,p)$
\cite{CG,EH2}
which is a nonstandard representation of the quantum group
$U_q(\widehat{sl_2})$.
\begin{eqnarray}
\[ R^{\mathrm{tr}}(q,p) \]_{ij}^{kl}
=p^{2(j-k)} \times
 \left\{
\begin{array}{cc}
q, & \mathrm{for} \ i=k \ge j=l, \\
q^{-1}, & \mathrm{for} \ i=k <  j=l, \\
-q+q^{-1}, & \mathrm{for} \ i<k<j, \ i+j=k+l, \\
q-q^{-1}, & \mathrm{for} \ j \le k<i, \
i+j=k+l, \\
0, & \mathrm{otherwise}.
\end{array}
\right.
\end{eqnarray}
This is nothing but the $R$-matrix in Definition \ref{definitionrkmatrixone}.
The CG-$R$-matrix $R^{\mathrm{tr}}(q,p)$ satisfies
the Yang-Baxter relation and the Hecke relation,
\begin{equation}
(\check{R}^{\mathrm{tr}}(q,p)-q)(\check{R}^{\mathrm{tr}}(q,p)+q^{-1})=0,
\end{equation}
thus can be used for representations for the generators
$T_1, \cdots, T_{n-1}$ as $T_j
=\check{R}^{\mathrm{tr}}_{j,j+1}(t,p), \ j=1,\cdots,n-1$.
The remaining step is to construct representations for the generators
$T_0$ and $T_{n}$.
This corresponds to finding solutions to the
reflection equations under the Hecke relation.
We achieve this by the following degeneration scheme \cite{EH}.
First, note there is a degeneration scheme to
the CG-$R$-matrix $R^{\mathrm{tr}}(q,p)$
from the Shibukawa-Ueno (SU) $R$-operator
$\mathcal{R}^{\mathrm{ell}}(\lambda)$,
which is an infinite-dimensional $R$-operator
acting on the function space.
This scheme can be summarized in the following diagram.
\begin{eqnarray}
\mathcal{R}^{\mathrm{ell}}(\lambda) 
\longrightarrow \mathcal{R}^{\mathrm{tr}}(\lambda)
\longrightarrow R^{\mathrm{tr}}(\lambda,q,p)
\longrightarrow R^{\mathrm{tr}}(q,p).
\end{eqnarray}
We start from the SU-$R$-operator
$\mathcal{R}^{\mathrm{ell}}(\lambda)$ \cite{SU}.
The first thing to do is to take the trigonometric limit of the $R$-operator.
Next, you twist the $R$-operator and take the finite-dimensional representation
in the trigonometric basis to get the spectral parameter dependent
CG-$R$-matrix $R^{\mathrm{tr}}(\lambda,q,p)$
from the infinite-dimensional $R$-operator.
Finally, one deBaxtetizes $R^{\mathrm{tr}}(\lambda,q,p)$
to obtain $R^{\mathrm{tr}}(q,p)$.
Note that the trigonometric $R$-matrix we consider here is different from
the standard $R$-matrix.
One should not reverse the order of the degeneration.
From the SU-$R$-operator, one obtains
the standard $R$-matrix by first taking the finite-dimensional representation
in the elliptic basis, and then taking the trigonometric limit,
but not by first taking the trigonometric limit of the SU-$R$-operator and then taking its finite-dimensional representation
in the trigonometric basis. What you get in this case is the
CG-$R$-matrix $R^{\mathrm{tr}}(\lambda,q,p)$.

Taking advantage of this fact of the degeneration and deBaxterization scheme
to obtain the $R$-matrix $R^{\mathrm{tr}}(\lambda,q,p)$,
we apply the same degeneration scheme to find solutions to the
reflection equations corresponding to the $R$-matrix 
$R^{\mathrm{tr}}(\lambda,q,p)$.
The scheme is given as follows.
\begin{eqnarray}
\mathcal{K}^{\mathrm{ell}}(\lambda) 
\longrightarrow \mathcal{K}^{\mathrm{tr}}(\lambda) 
\longrightarrow K^{\mathrm{tr}}(\lambda,r,s)
\longrightarrow K^{\mathrm{tr}}(r,s).
\end{eqnarray}
We start from the elliptic Hikami-Komori (HK) $K$-operator
$\mathcal{K}^{\mathrm{ell}}(\lambda)$ \cite{Hi,HK},
which corresponds to the solution to the reflection equation
of the SU-$R$-operator $\mathcal{R}^{\mathrm{ell}}(\lambda)$.
First, we take the trigonometric limit
$\mathcal{K}^{\mathrm{tr}}(\lambda)$ of the
elliptic $K$-operator $\mathcal{K}^{\mathrm{ell}}(\lambda)$.
Next, we twist the $K$-operator and take the finite-dimensional representation
in the trigonometric basis to obtain
$K^{\mathrm{tr}}(\lambda,r,s)$.
Finally, we deBaxterize $K^{\mathrm{tr}}(\lambda,r,s)$
to get $K^{\mathrm{tr}}(r,s)$,
whose explicit matrix elements are given by
\begin{eqnarray}
\[ K^{\mathrm{tr}}(r,s) \]_{j}^{k}
=s^{j-k} \times
 \left\{
\begin{array}{cc}
-r^{-1}, & \mathrm{for} \ j \le k, \ j+k=N-1, \\
-r, & \mathrm{for} \ j>k, \ j+k=N-1, \\
r-r^{-1}, & \mathrm{for} \ j \le k<N-1-j, \\
-r+r^{-1}, & \mathrm{for} \ N-1-j<k<j, \\
0, & \mathrm{otherwise}.
\end{array}
\right.
\end{eqnarray}
which is the matrix given in Definition \ref{definitionrkmatrixone}.
One can show that this $K$-matrix satisfies the Hecke relation
\begin{equation}
(K^{\mathrm{tr}}(r,s)-r)(K^{\mathrm{tr}}(r,s)+r^{-1})=0,
\end{equation}
and thus can be served as representations for $T_0$ and $T_n$
as $T_0=K_1^{\mathrm{tr}}(t_0,s_0)$ and 
$T_n=K_n^{\mathrm{tr}}(t_n,s_n)$.
One can easily see the remaining relation
\eqref{relseven} holds by looking at which spaces
the $R$-matrices and the $K$-matrices act on nontrivially.
This ends the proof that the matrices given in
Definition \ref{definitionrkmatrixone} satisfy all the
defining relations for the affine Hecke algebra of type $C$,
and one can use them as representations for the generators,
thus proving Theorem \ref{thone}.
The details of the calculation are given in the next section.
\subsection{Rational representation}
For the special case $t=t_0=t_1=1$ of type $C$ affine Hecke algebra 
$H_n(t,t_0,t_1)$,
one can obtain another representation for the algebra
from the rational Jordanian $R$-matrix and its corresponding $K$-matrix.
Since the degeneration \cite{EH} and deBaxterization scheme is similar,
let us point out the differences.
On the deBaxterization scheme, one replaces
the decompostions of the $R$-matrix $R(\lambda)$ \eqref{baxterone}
and $K$-matrix $K(\lambda)$ \eqref{baxterthree} by
\begin{align}
&R(\lambda)=f(\lambda)(P+g\lambda R), \\
&K(\lambda)=a(\lambda)(I+b \lambda K),
\end{align}
and the Hecke relations in \eqref{baxtertwo} and \eqref{baxterfour} by
\begin{align}
\check{R}^2-1=0, \\
K^2-1=0.
\end{align}
The diagram to get the rational Jordanian $R$-matrix
$R^{\mathrm{ra}}(\kappa,h)$ from the SU-$R$-operator
$\mathcal{R}^{\mathrm{ell}}(\lambda)$ is changed as
\begin{eqnarray}
\mathcal{R}^{\mathrm{ell}}(\lambda)
\longrightarrow \mathcal{R}^{\mathrm{tr}}(\lambda)
\longrightarrow \mathcal{R}^{\mathrm{ra}}(\lambda)
\longrightarrow R^{\mathrm{ra}}(\lambda,\kappa,h)
\longrightarrow R^{\mathrm{ra}}(\kappa,h).
\end{eqnarray}
The difference from the procedure to obtain the
trigonometric CG-$R$-matrix is that
we degenerate furthermore the infinite-dimensional $R$-operator
from the trigonometric one $\mathcal{R}^{\mathrm{tr}}(\lambda)$ to
the rational one $\mathcal{R}^{\mathrm{ra}}(\lambda)$.
Then we take the finite-dimensional representation of
the $R$-operator to get $R^{\mathrm{ra}}(\lambda,\kappa,h)$,
and deBaxterize it to obtain $R^{\mathrm{ra}}(\kappa,h)$.

Correspondingly, the $K$-matrix for the reflection equation
of the Jordanian $R$-matrix can be obtained in the following procedure:
\begin{eqnarray}
\mathcal{K}^{\mathrm{ell}}(\lambda)
\longrightarrow \mathcal{K}^{\mathrm{tr}}(\lambda) 
\longrightarrow \mathcal{K}^{\mathrm{ra}}(\lambda) 
\longrightarrow K^{\mathrm{ra}}(\lambda,\nu,g)
\longrightarrow K^{\mathrm{ra}}(\nu,g).
\end{eqnarray}
First, we degenerate the trigonometric HK-$K$-operator
futhermore to the rational $K$-operator.
Then we take the finite-dimensional representation of the rational $K$-matrix
$K^{\mathrm{ra}}(\lambda,\nu,g)$, and deBaxterize
to obtain the rational $K$-matrix without spectral parameter
$K^{\mathrm{ra}}(\nu,g)$.
The explicit forms of the rational Jordanian $R$-matrix and $K$-matrix
are given in Definition \ref{definitionrkmatrixtwo},
and can be served as a representation for $H_n(1,1,1)$.
The details are given in the section 5.
\section{Trigonometric representation}
In this and the next sections, we give details of the proof
outlined in the last section to construct representations.
In this section,
we consider the nonstandard trigonometric representations
of type $C$ affine Hecke algebra.
The term trigonometric comes from the fact that the representation
comes from quantum integrable models of trigonometric type.
For completeness, we first review the
degeneration scheme \cite{EH} from the Shibukawa-Ueno (SU)
elliptic $R$-operator \cite{SU}
to the Cremmer-Gervais (CG) trigonometric $R$-matrix.
Then we give the details of obtaining solutions
to the reflection equation from the Hikami-Komori (HK) $K$-operator.
We also compare the obtained solution
with our former result on the full solution
for $N=3$ \cite{MY}.
\subsection{Cremmer-Gervais $R$-matrix from Shibukawa-Ueno $R$-operator}
In this section, we review how the Cremmer-Gervais $R$-matrix is
extracted from the Shibukawa-Ueno $R$-operator.
\begin{definition} \cite{SU}
Let $\mathcal{M}$ be a space of meromorphic functions of $\zeta$ on
$\mathbf{C}^n$. The Shibukawa-Ueno $R$-operator $\mathcal{R}^{\mathrm{ell}}_{jk}(\lambda)
\in \mathrm{End}(\mathcal{M})$ is defined as
\begin{equation}
\mathcal{R}^{\mathrm{ell}}_{jk}(\lambda)=\sigma_\lambda(\zeta_j-\zeta_k;\tau)s_{jk}-
\sigma_\kappa(\zeta_j-\zeta_k;\tau),
\end{equation}
where
\begin{align}
&\sigma_\nu(\zeta;\tau)=\frac{\theta_1(\zeta+\nu;\tau)\theta_1^\prime(0;\tau)}
{\theta_1(\zeta;\tau)\theta_1(\nu;\tau)}, \ \theta_1(\zeta;\tau)
=-\theta
\left[
\begin{array}{c}
\frac{1}{2} \\
\frac{1}{2}
\end{array}
\right] (\zeta;\tau), \\
&\theta
\left[
\begin{array}{c}
a \\
b
\end{array}
\right] (\zeta;\tau)
=\sum_{m \in \mathbf{Z}} \mathrm{exp}
\Big(
\pi \mathrm{i} \tau (m+a)^2
+2 \pi \mathrm{i} (m+a)(\zeta+b)
\Big), \\
&(s_{jk}f)(\zeta,\cdots,\zeta_j,\cdots,\zeta_k,\cdots,\zeta_n)=f(\zeta_1,\cdots,\zeta_k,\cdots,\zeta_j,\cdots,\zeta_n).
\end{align}
\hspace*{\fill}    \sq
\end{definition}
The SU-$R$-operator
can be twisted \cite{FP}
using the shift operator
\begin{equation}
(T_{j}(\mu)f)(\zeta_1,\cdots,\zeta_j,\cdots,\zeta_n)
=f(\zeta_1,\cdots,\zeta_j+\mu,\cdots,\zeta_n),
\end{equation}
as
\begin{equation}
\widetilde{\mathcal{R}}^{\mathrm{ell}}_{jk}(\lambda):=T_j(-\beta)T_k(\beta)\mathcal{R}_{jk}^{\mathrm{ell}}
(\lambda)
T_j(-\beta)T_k(\beta).
\end{equation}
\begin{theorem} \cite{SU}
The Shibukawa-Ueno $R$-operator $\widetilde{\mathcal{R}}_{jk}^{\mathrm{ell}}(\lambda)$
satisfies the Yang-Baxter relation
\begin{equation}
\widetilde{\mathcal{R}}_{12}(\lambda)
\widetilde{\mathcal{R}}_{13}(\lambda+\mu)
\widetilde{\mathcal{R}}_{23}(\mu)
=\widetilde{\mathcal{R}}_{23}(\lambda)
\widetilde{\mathcal{R}}_{13}(\lambda+\mu)
\widetilde{\mathcal{R}}_{12}(\mu).
\end{equation}
\hspace*{\fill}    \sq
\end{theorem}
The action of $\widetilde{\mathcal{R}}_{12}^{\mathrm{ell}}(\lambda)$
on $f(\zeta_1,\zeta_2)$ is explicitly given as
\begin{align}
\widetilde{\mathcal{R}}_{12}^{\mathrm{ell}}(\lambda)f(\zeta_1,\zeta_2)
=&\sigma_\lambda(\zeta_1-\zeta_2-2\beta;\tau)f(\zeta_2,\zeta_1)
-\sigma_\kappa(\zeta_1-\zeta_2;\tau)f(\zeta_1-2 \beta,\zeta_2+2 \beta).
\end{align}
The CG-$R$-matrix
can be obtained from the SU-$R$-operator as follows.
First, one takes the trigonometric limit of the elliptic $R$-operator
\begin{equation}
\widetilde{\mathcal{R}}_{12}^{\mathrm{tr}}(\lambda)
:=
(2 \pi \mathrm{i})^{-1}
\mathrm{lim}_{\tau \to \infty}\widetilde{\mathcal{R}}_{12}^{\mathrm{ell}}(\lambda).
\end{equation}
The action of the trigonometric $R$-operator on the function 
$f(\zeta_1,\zeta_2)$ is given by
\begin{align}
\widetilde{\mathcal{R}}_{12}^{\mathrm{tr}}(\lambda)f(\zeta_1,\zeta_2)
=&\frac{zw_1-z^{-1}p^2 w_2}{(z-z^{-1})(w_1-p^2 w_2)}f(\zeta_2,\zeta_1)
-\frac{qw_1-q^{-1}p^2 w_2}{(q-q^{-1})(w_1-p^2 w_2)}
f(\zeta_1-2 \beta,\zeta_2+2 \beta),
\end{align}
where we have defined
$w_j=\mathrm{e}^{2 \pi \mathrm{i} \zeta_j}, \ 
z=\mathrm{e}^{\pi \mathrm{i} \lambda}, \
p=\mathrm{e}^{2 \pi \mathrm{i} \beta}, \
q=\mathrm{e}^{\pi \mathrm{i} \kappa}
$.
Next, one takes the finite-dimensional representation.
We restrict the space of functions $\mathcal{M}$ to the finite-dimensional
subspace $V_N^{\mathrm{tr}} \otimes V_N^{\mathrm{tr}}$
where
\begin{equation}
V_N^{\mathrm{tr}}=\mathrm{Span}
\Big\{
\phi_k(\zeta)=\mathrm{e}^{\pi \mathrm{i} (2k-N+1) \zeta} \ \Big| \
k=0,1,\cdots,N-1 
\Big\}.
\end{equation}
Calculating the matrix elements
of the trigonometric $R$-operator explicitly
using $\phi_k(\zeta), \ k=0,1,\cdots,N-1$ as the basis,
one gets \cite{EH}
\begin{align}
\widetilde{\mathcal{R}}_{12}^{\mathrm{tr}}(\lambda) 
\phi_i(\zeta_1) \phi_j(\zeta_2)
=\sum_{k,l=0}^{N-1}
\[ R^{\mathrm{tr}}(\lambda,q,p) \]_{ij}^{kl} \phi_k(\zeta_1) \phi_l(\zeta_2),
\end{align}
with \\
$ \[ R^{\mathrm{tr}}(\lambda,q,p) \]_{ij}^{kl}  $
\begin{eqnarray}
=p^{2(j-k)} \times
 \left\{
\begin{array}{cc}
(q z^{-1}-q^{-1} z) \slash (q-q^{-1})(z-z^{-1}), & \mathrm{for} \ i=j=k=l, \\
-q^{\mathrm{sgn}(i-j)} \slash (q-q^{-1}), & \mathrm{for} \ i=k \neq  j=l, \\
z^{\mathrm{sgn}(j-i)} \slash (z-z^{-1}), & \mathrm{for} \ l=i \neq  k=j, \\
\mathrm{sgn}(j-i), & \mathrm{for} \ \mathrm{min}(i,j) < k < \mathrm{max}(i,j), \
i+j=k+l, \\
0, & \mathrm{otherwise}. \label{spectral-cremmer-gervais}
\end{array}
\right.
\end{eqnarray}
Note that $z$ is defined as $z=\mathrm{e}^{\pi \mathrm{i} \lambda}$.
The $R$-matrix $R^{\mathrm{tr}}(\lambda,q,p)$
is called the Cremmer-Gervais $R$-matrix \cite{CG,EH}.
From its construction from the SU-$R$-operator,
the CG-$R$-matrix satisfies the Yang-Baxter relation
\begin{equation}
R_{12}(\lambda)R_{13}(\lambda+\mu)R_{23}(\mu)=
R_{23}(\mu)R_{13}(\lambda+\mu)R_{12}(\lambda). \label{yang-baxter-ver2}
\end{equation}
We now apply the following lemma to extract the constant $R$-matrix.
\begin{lemma} \label{gettingridoflemma}
Let an $R$-matrix $R(\lambda)$ satisfying the Yang-Baxter relation
\eqref{yang-baxter-ver2}
can be decomposed as
\begin{equation}
R(\lambda)=f(\lambda)(P+g(\lambda)R), \ g(\lambda)=
\frac{\mathrm{e}^{-2 \pi \mathrm{i} \lambda}-1}{q-q^{-1}}, \ f(\lambda) \not\equiv 0,
\label{decomposition-ver2}
\end{equation}
and the $\lambda$-independent $R$ satisfies the Hecke relation
\begin{equation}
(\check{R}-q)(\check{R}+q^{-1})=0, \label{hecke-ver2}
\end{equation}
where $\check{R}=RP$.
The $\lambda$-independent $R$ satisfies the Yang-Baxter relation without spectral parameter.
\begin{equation}
R_{12}R_{13}R_{23}=
R_{23}R_{13}R_{12}. \label{ybr-ver2}
\end{equation}
\hspace*{\fill}    \sq
\end{lemma}
We now apply this lemma.
One finds the $R$-matrix $R^{\mathrm{tr}}(\lambda,q,p)$
\eqref{spectral-cremmer-gervais} can be decomposed as
\begin{equation}
R^{\mathrm{tr}}(\lambda,q,p)
=
\frac{1}{1-\mathrm{e}^{-2 \pi \mathrm{i} \lambda}}
\Bigg( P+
\frac{\mathrm{e}^{-2 \pi \mathrm{i} \lambda}-1}{q-q^{-1}}
R^{\mathrm{tr}}(q,p) \Bigg),
\label{explicitdecomposition}
\end{equation}
with the $\lambda$-independent $R$-matrix given by
\begin{eqnarray}
\[ R^{\mathrm{tr}}(q,p) \]_{ij}^{kl}
=p^{2(j-k)} \times
 \left\{
\begin{array}{cc}
q, & \mathrm{for} \ i=k \ge j=l, \\
q^{-1}, & \mathrm{for} \ i=k <  j=l, \\
-q+q^{-1}, & \mathrm{for} \ i<k<j, \ i+j=k+l, \\
q-q^{-1}, & \mathrm{for} \ j \le k<i, \
i+j=k+l, \\
0, & \mathrm{otherwise}.
\end{array}
\right.
\end{eqnarray}
Since this $R$-matrix can be shown to satisfy the Hecke relation
\eqref{hecke-ver2}, we can apply the above lemma and find that the 
$R$-matrix $R^{\mathrm{tr}}(q,p)$ satisfy the Yang-Baxter relation
\eqref{ybr-ver2} which is equivalent to the braid relation
\begin{equation}
\check{R}_{12}\check{R}_{23}\check{R}_{12}
=\check{R}_{23}\check{R}_{12}\check{R}_{23}. \label{braid-ver2}
\end{equation}
This braid relation \eqref{braid-ver2} together with the Hecke relation
\eqref{hecke-ver2} can be identified with the defining relations
\eqref{relfour} and \eqref{reltwo} of type $C$ affine Hecke algebra,
hence we can use the CG-$R$-matrix
multiplied by the permutation matrix
$\check{R}^{\mathrm{tr}}(q,p)$
as a representation for the generators
$T_1,\cdots,T_{n-1}$ of the affine Hecke algebra of type $C$:
\begin{equation}
T_j=\check{R}_j^{\mathrm{tr}}(t,p), \ j=1,\cdots,n-1. \label{forthone1}
\end{equation}

\subsection{$K$-matrix from Hikami-Komori $K$-operator}
We now apply the same degenataion scheme to find solutions
to the reflection equation of the $R$-matrix $R^{\mathrm{tr}}(q,p)$.
We start from the elliptic Hikami-Komori (HK) $K$-operator.
\begin{definition} \cite{Hi,HK}
Let $\mathcal{M}$ be a space of meromorphic functions of $\zeta$ on
$\mathbf{C}^n$. The Hikami-Komori $K$-operator $\mathcal{K}_{j}^{\mathrm{ell}}(\lambda)
\in \mathrm{End}(\mathcal{M})$ is defined as
\begin{equation}
\mathcal{K}_{j}^{\mathrm{ell}}(\lambda)=\sigma_\nu(\zeta;\tau)s_{j}-
\sigma_{2 \lambda}(\zeta;\tau),
\end{equation}
where
\begin{align}
(s_{j}f)(\zeta,\cdots,\zeta_j,\cdots,\zeta_n)=f(\zeta_1,\cdots,-\zeta_j,\cdots,\zeta_n).
\end{align}
\hspace*{\fill}    \sq
\end{definition}
The HK-$K$-operator
can be twisted using the shift operator as
\begin{equation}
\widetilde{\mathcal{K}}_{j}^{\mathrm{ell}}(\lambda):=T_j(-\gamma)
\mathcal{K}_{j}^{\mathrm{ell}}(\lambda) T_j(\gamma).
\end{equation}
\begin{theorem} \cite{Hi,HK}
The Hikami-Komori $K$-operator $\widetilde{\mathcal{K}}_{j}^{\mathrm{ell}}(\lambda)$
is a solution to the reflection equation of the Shibukawa-Ueno
$R$-operator
\begin{equation}
\widetilde{\mathcal{R}}_{12}(\lambda-\mu)
\widetilde{\mathcal{K}}_{1}(\lambda)
\widetilde{\mathcal{R}}_{21}(\lambda+\mu)
\widetilde{\mathcal{K}}_{2}(\mu)
=
\widetilde{\mathcal{K}}_{2}(\mu)
\widetilde{\mathcal{R}}_{12}(\lambda+\mu)
\widetilde{\mathcal{K}}_{1}(\lambda)
\widetilde{\mathcal{R}}_{21}(\lambda-\mu).
\end{equation}
\hspace*{\fill}    \sq
\end{theorem}
The action of $\widetilde{\mathcal{K}}^{\mathrm{ell}}(\lambda)$
on $f(\zeta)$ is explicitly given as
\begin{align}
\widetilde{\mathcal{K}}(\lambda)f(\zeta)
=&\sigma_\nu(\zeta-\gamma;\tau)f(-\zeta+2 \gamma)
-\sigma_{2 \lambda}(\zeta-\gamma;\tau)f(\zeta).
\end{align}
Now we calculate the $K$-matrix corresponding to the
CG-$R$-matrix
from the HK-$K$-operator, following the same line
as the previous subsection.
First, one takes the trigonometric limit of the elliptic $K$-operator
\begin{equation}
\widetilde{\mathcal{K}}^{\mathrm{tr}}(\lambda)
:=
(2 \pi \mathrm{i})^{-1}
\mathrm{lim}_{\tau \to \infty}\widetilde{\mathcal{K}}^{\mathrm{ell}}(\lambda).
\end{equation}
The action of the untwisted ($\gamma=0$) trigonometric $K$-operator
on the function 
$f(\zeta)$ is given by
\begin{align}
\mathcal{K}^{\mathrm{tr}}(\lambda)f(\zeta)
=&\frac{r-wr^{-1}}{(r-r^{-1})(w-1)}f(-\zeta)
-\frac{w z^2-z^{-2}}{(z^2-z^{-2})(w-1)}
f(\zeta),
\end{align}
where we have defined
$w=\mathrm{e}^{2 \pi \mathrm{i} \zeta}, \ 
z=\mathrm{e}^{\pi \mathrm{i} \lambda}, \
r=\mathrm{e}^{- \pi \mathrm{i} \nu}
$.
Next, one takes the finite-dimensional representation.
We restrict the space of functions $\mathcal{M}$ to the finite-dimensional
subspace $V_N^{\mathrm{tr}}$.
Calculating explicitly the matrix elements
of the trigonometric $K$-operator
using $\phi_k(\zeta), \ k=0,1,\cdots,N-1$ as the basis,
one gets the following:
\begin{proposition}
The matrix elements
$\[ K^{\mathrm{tr}}(\lambda,r,s) \]_{j}^{k}$
of the trigonometric $K$-operator
$\widetilde{\mathcal{K}}^{\mathrm{tr}}(\lambda)$
in the basis $\phi_k(\zeta), \ k=0,1,\cdots,N-1$
\begin{align}
\widetilde{\mathcal{K}}^{\mathrm{tr}}(\lambda) 
\phi_j(\zeta)
=\sum_{k=0}^{N-1}
\[ K^{\mathrm{tr}}(\lambda,r,s) \]_{j}^{k} \phi_k(\zeta),
\end{align}
is given by \\
$ \[ K^{\mathrm{tr}}(\lambda,r,s) \]_{j}^{k}  $
\begin{eqnarray}
=s^{j-k} \times
 \left\{
\begin{array}{cc}
(r^{-1} z^{-2}-r z^2) \slash (r-r^{-1})(z^2-z^{-2}), & \mathrm{for} \ j=k=(N-1)/2, \\
-z^{2 \mathrm{sgn}(2j-N+1)} \slash (z^2-z^{-2}), & \mathrm{for} \ j=k \neq (N-1)/2, \\
-r^{\mathrm{sgn}(2j-N+1)} \slash (r-r^{-1}), & \mathrm{for} \ k=N-1-j, \  j \neq (N-1)/2, \\
\mathrm{sgn}(N-1-2j), & \mathrm{for} \ \mathrm{min}(j,N-1-j) < k < \mathrm{max}(j,N-1-j), \\
0, & \mathrm{otherwise}. \label{spectral-cremmer-gervais-k}
\end{array}
\right.
\end{eqnarray}
\hspace*{\fill}    \sq
\end{proposition}
{\it Proof}. \\
We first consider the case with no twist $\gamma=0$.
The case with general $\gamma$ can be obtained from $\gamma=0$
in a simple way.

We act $\mathcal{K}^{\mathrm{tr}}(\lambda)$ on $\phi_j(\zeta)$
\begin{align}
\mathcal{K}^{\mathrm{tr}}(\lambda) \phi_j(\zeta)
=&\frac{r-wr^{-1}}{(r-r^{-1})(w-1)}
\mathrm{e}^{\pi \mathrm{i} (N-2j-1) \zeta}
-\frac{w z^2-z^{-2}}{(z^2-z^{-2})(w-1)}
\mathrm{e}^{\pi \mathrm{i} (2j-N+1) \zeta} \nonumber \\
=&
\mathrm{e}^{\pi \mathrm{i} (1-N) \zeta}
\Bigg\{
\frac{r-wr^{-1}}{(r-r^{-1})(w-1)}
w^{N-1-j}
-\frac{w z^2-z^{-2}}{(z^2-z^{-2})(w-1)}
w^{j}
\Bigg\} \nonumber \\
=&\frac{\mathrm{e}^{\pi \mathrm{i} (1-N) \zeta}}{(r-r^{-1})(z^2-z^{-2})}
\times
\frac{(z^2-z^{-2})(r-wr^{-1})w^{N-1-j}-(r-r^{-1})(wz^2-z^{-2})w^j}{w-1}.
\label{int1}
\end{align}
We reorganize the second factor
in the last line into polynomials in $w$ as follows.
\begin{align}
&\frac{(z^2-z^{-2})(r-wr^{-1})w^{N-1-j}-(r-r^{-1})(wz^2-z^{-2})w^j}{w-1}
\nonumber \\
=&\frac{1}{w-1}
\left|
\begin{array}{cc}
w^{N-1-j} & w^j \\
(r-r^{-1})(w z^2-z^{-2}) & (z^2-z^{-2})(r-wr^{-1})
\end{array}
\right| \nonumber \\
=&\frac{1}{w-1}
\left|
\begin{array}{cc}
w^{N-1-j}-w^j & w^j \\
(w-1)(r z^2-r^{-1} z^{-2}) & (z^2-z^{-2})(r-wr^{-1})
\end{array}
\right| \nonumber \\
=&(z^2-z^{-2})(r-wr^{-1}) \frac{w^{N-1-j}-w^j}{w-1}
-w^j (r z^2-r^{-1} z^{-2}) \nonumber \\
=&(z^2-z^{-2})(r-wr^{-1}) \sum_l \epsilon(j,N-1-j,l)w^l
-w^j (r z^2-r^{-1} z^{-2}) \nonumber \\
=&(z^2-z^{-2}) \sum_l 
\{
r \epsilon(j,N-1-j,l)w^l
-r^{-1} \epsilon(j,N-1-j,l)w^{l+1}
\}
-w^j (r z^2-r^{-1} z^{-2}) \nonumber \\
=&(z^2-z^{-2}) \sum_l
\{
 r \epsilon(j,N-1-j,l)-r^{-1} \epsilon(j,N-1-j,l-1)
\}
w^l
-w^j (r z^2-r^{-1} z^{-2}). \label{int2}
\end{align}
The sum in the last line can be explicitly calculated using
the definition of $\epsilon(i,j,k)$
\begin{eqnarray}
\epsilon(i,j,k)
=\left\{
\begin{array}{cc}
1, & \mathrm{for} \ i \le k < j, \\
-1, & \mathrm{for} \ j \le k < i, \\
0, & \mathrm{otherwise},
\end{array}
\right.
\end{eqnarray}
as
\begin{align}
&r \epsilon(j,N-1-j,l)-r^{-1} \epsilon(j,N-1-j,l-1) \nonumber \\
=& \left\{
\begin{array}{cc}
(r-r^{-1}) \mathrm{sgn}(N-1-2j), & \mathrm{for} \ \mathrm{min}(j,N-1-j)<l<
\mathrm{max}(j,N-1-j), \\
r^{\mathrm{sgn}(N-1-2j)}, & \mathrm{for} \ l=j \neq N-1-j, \\
-r^{\mathrm{sgn}(2j-N+1)}, & \mathrm{for} \ l=N-1-j \neq j, \\
0, & \mathrm{otherwise}. \label{int3}
\end{array}
\right.
\end{align}
Combining \eqref{int1}, \eqref{int2} and \eqref{int3}
gives the proof of the proposition for the case with no twist $\gamma=0$.
The case with nonzero twist $\gamma$ can be included through a simple relation.
Since the action of the shift operator on the trigonometric basis is diagonal
\begin{align}
T(\gamma) \phi_j(\zeta)=\mathrm{e}^{\pi \mathrm{i} \tau (2k-N+1)}
\phi_j(\zeta),
\end{align}
one has the following simple relation for the matrix elements
between the twisted and nontwisted $K$-matrices
$
\widetilde{\mathcal{K}}^{\mathrm{tr}}(\lambda)
=
T(-\gamma)
\mathcal{K}^{\mathrm{tr}}(\lambda)
T(\gamma)
$
\begin{align}
[K^{\mathrm{tr}}(\lambda,r,s)]_j^k
=s^{j-k} [K^{\mathrm{tr}}(\lambda,r)]_j^k,
\end{align}
where $s=\mathrm{e}^{2 \pi \mathrm{i} \gamma}$, which concludes the proof
including twist.
\hspace*{\fill}    \wsq

\begin{proposition}
The matrix $K^{\mathrm{tr}}(\lambda,r,s)$
is a solution to the reflection equation
of the trigonometric CG-$R$-matrix $R^{\mathrm{tr}}(\lambda,r,s)$.
\hspace*{\fill}    \sq
\end{proposition}
{\it Proof}. \\
This follows from the fact that the $K$-matrix
is constructed as a degeneration from the HK-$K$-operator,
which is a solution to the reflection equation of the
SU-$R$-operator. \hspace*{\fill}    \wsq \\
So far, we found a solution $K^{\mathrm{tr}}(\lambda,r,s)$
to the reflection equation of the CG-$R$-matrix
$R^{\mathrm{tr}}(\lambda,q,p)$ with spectral parameter
\begin{equation}
R_{12}(\lambda-\mu)
K_{1}(\lambda)
R_{21}(\lambda+\mu)
K_{2}(\mu)
=
K_{2}(\mu)
R_{12}(\lambda+\mu)
K_{1}(\lambda)
R_{21}(\lambda-\mu). \label{reflectionforlemma}
\end{equation}
To extract representations for the generators $T_0$ and $T_n$,
we use the following lemma.
\begin{lemma} \label{auxlemmatrigonometric}
Let $R(\lambda)$ be an $R$-matrix which satisfies the properties in
Lemma \ref{gettingridoflemma}.
If the corresponding $K$-matrix $K(\lambda)$ satisfying
the reflection equation
\eqref{reflectionforlemma}
can be decomposed as
\begin{equation}
K(\lambda)=a(\lambda)(I+b(\lambda)K), \ b(\lambda)
=\frac{\mathrm{e}^{-4 \pi \mathrm{i} \lambda}-1}{r-r^{-1}}, \ a(\lambda) \not\equiv 0,
\label{decompositiontrigonometric}
\end{equation}
and the $\lambda$-independent $K$ satisfies the Hecke relation
\begin{equation}
(K-r)(K+r^{-1})=0, \label{hecketrigonometric}
\end{equation}
then the $\lambda$-independent $K$ satisfies the reflection equations
without spectral parameter
\begin{align}
&\check{R}_{12}K_{1}\check{R}_{12}K_{1}
=K_{1}\check{R}_{12}K_{1}\check{R}_{12}, \label{whatweneed1} \\
&\check{R}_{12}K_{2}\check{R}_{12}K_{2}
=K_{2}\check{R}_{12}K_{2}\check{R}_{12}. \label{whatweneed2}
\end{align}
\hspace*{\fill}    \sq
\end{lemma}
{\it Proof}. \\
Multiplying both sides of the
spectral parameter dependent reflection equation
\eqref{reflectionforlemma} by the permutation operators $P$ and
inserting the decomposition relations \eqref{decomposition-ver2},
\eqref{decompositiontrigonometric}
and Hecke relations \eqref{hecke-ver2},
\eqref{hecketrigonometric}
into it,
one finds the coefficients of the terms $K \check{R}$
and $\check{R} K$ have the following form:
\begin{align}
&(b(\mu)+b(\lambda))g(\lambda-\mu)+(b(\mu)-b(\lambda))g(\lambda+\mu)
\nonumber \\
&+(q-q^{-1})b(\mu)g(\lambda+\mu)g(\lambda-\mu)
+(r-r^{-1})b(\lambda)b(\mu)g(\lambda-\mu).
\end{align}
We can show by explicit calculation that this becomes zero.
Cancelling out these vanishing terms,
The remaining equations are nothing but
the reflection equations without spectral parameters
\eqref{whatweneed1}, \eqref{whatweneed2}. \hspace*{\fill}    \wsq \\
We now apply this lemma to extract the constant $K$-matrix.

\begin{proposition} \label{auxproptrigonometric}
The $K$-matrix $K^{\mathrm{tr}}(r,s)$
whose matrix elements are explicitly given by
\begin{eqnarray}
\[ K^{\mathrm{tr}}(r,s) \]_{j}^{k}
=s^{j-k} \times
 \left\{
\begin{array}{cc}
-r^{-1}, & \mathrm{for} \ j \le k, \ j+k=N-1, \\
-r, & \mathrm{for} \ j>k, \ j+k=N-1, \\
r-r^{-1}, & \mathrm{for} \ j \le k<N-1-j, \\
-r+r^{-1}, & \mathrm{for} \ N-1-j<k<j, \\
0, & \mathrm{otherwise},
\end{array}
\right. \label{forprop}
\end{eqnarray}
is a solution to the constant reflection equation
of the CG-$R$-matrix $R^{\mathrm{tr}}(p,q)$.
\hspace*{\fill}    \sq
\end{proposition}
{\it Proof}. \\
The $K$-matrix $K^{\mathrm{tr}}(\lambda,r,s)$
satisfying the reflection equation \eqref{reflectionforlemma}
can be decomposed as
\begin{align}
K^{\mathrm{tr}}(\lambda,r,s)=\mathrm{e}^{2 \pi \mathrm{i} \lambda}
(r^{-1}-r) 
\Bigg(
I+
\frac{\mathrm{e}^{-4 \pi \mathrm{i} \lambda}-1}{r-r^{-1}}K^{\mathrm{tr}}(r,s) \Bigg),
\label{propdecomposition}
\end{align}
with the $\lambda$-dependent $K$-matrix
$K^{\mathrm{tr}}(r,s)$ explicitly given by \eqref{forprop}.
To apply Lemma \ref{auxlemmatrigonometric},
one also needs to show that $K^{\mathrm{tr}}(r,s)$ satisfies the Hecke relation
\eqref{hecketrigonometric}.
This follows by comparing the expression
\begin{align}
K^{\mathrm{tr}}(\lambda,r,s)K^{\mathrm{tr}}(-\lambda,r,s)
=\{ (r-r^{-1})^2-(\mathrm{e}^{2 \pi \mathrm{i} \lambda}-
\mathrm{e}^{-2 \pi \mathrm{i} \lambda})^2 \} I,
\end{align}
which can be calculated using the operator expression for
$K^{\mathrm{tr}}(\lambda,r,s)$, and comparing with another expression
\begin{align}
K^{\mathrm{tr}}(\lambda,r,s)K^{\mathrm{tr}}(-\lambda,r,s)
=(r-r^{-1})^2 I+
(\mathrm{e}^{2 \pi \mathrm{i} \lambda}
-\mathrm{e}^{-2 \pi \mathrm{i} \lambda})^2
((r-r^{-1})K-K^2),
\end{align}
obtained from the decomposition
\eqref{propdecomposition}. \\
We have shown all the conditions the $K$-matrices must satisfy
to apply Lemma \ref{auxlemmatrigonometric},
hence the proposition follows from the lemma.
\hspace*{\fill}    \wsq \\
Lemma \ref{auxlemmatrigonometric} and Proposition \ref{auxproptrigonometric}
shows that the $K$-matrix
$K^{\mathrm{tr}}(\lambda,r,s)$ satisfies
the reflection equations \eqref{whatweneed1}, \eqref{whatweneed2}
and the Hecke relation \eqref{hecketrigonometric}
which can be identified with the
defining relations \eqref{relfour}, \eqref{relsix},
\eqref{relone} and \eqref{relthree}
of type $C$ affine Hecke algebra,
one obtains a representation for $T_0$ and $T_n$
in terms of the $K$-matrix:
\begin{align}
T_0=K^{\mathrm{tr}}_1(t_0,s_0), \label{forthone2} \\
T_n=K^{\mathrm{tr}}_n(t_n,s_n). \label{forthone3}
\end{align}
The representation for the generators in terms
of the constant CG-$R$-matrix and
its assocaited $K$-matrix
\eqref{forthone1}, \eqref{forthone2} and \eqref{forthone3}
proves Theorem \ref{thone}.
\subsection{Another representation}
One can obtain another representation starting from
another elliptic $K$-operator \cite{HK}
\begin{equation}
\overline{\mathcal{K}}_{j}(\lambda)
=\overline{\sigma}_\nu(\zeta;\tau)s_{j}-
\overline{\sigma}_{2 \lambda}(\zeta;\tau),
\end{equation}
where
\begin{equation}
\overline{\sigma}_\nu(\zeta;\tau)
=\frac{\theta_2(\zeta+\nu;\tau)\theta_1^\prime(0;\tau)}
{\theta_2(\zeta;\tau)\theta_1(\nu;\tau)}, \
\theta_2(\zeta;\tau)=
\theta
\left[
\begin{array}{c}
\frac{1}{2} \\
0
\end{array}
\right] (\zeta;\tau),
\end{equation}
which is another solution satisfying the reflection equation of the
SU-$R$-operator.
We just present the results for the explicit matrix elements
of the finite-dimensional $K$-matrix for $N$ odd which can be obtained
in the same line as the previous subsection.
\begin{proposition}
The following $K$-matrices
$\overline{K}^{\mathrm{tr}}(\lambda,r,s)$
and
$\overline{K}^{\mathrm{tr}}(r,s)$
is a solution to the reflection equation
of the CG-$R$-matrix with and without spectral parameter,
respectively. \\
$ \[ \overline{K}^{\mathrm{tr}}(\lambda,r,s) \]_{j}^{k}  $
\begin{eqnarray}
=s^{j-k} \times
 \left\{
\begin{array}{cc}
(r^{-1} z^{-2}-r z^2) \slash (r-r^{-1})(z^2-z^{-2}), & \mathrm{for} \ j=k=(N-1)/2, \\
-z^{2 \mathrm{sgn}(2j-N+1)} \slash (z^2-z^{-2}), & \mathrm{for} \ j=k \neq (N-1)/2, \\
-r^{\mathrm{sgn}(2j-N+1)} \slash (r-r^{-1}), & \mathrm{for} \ k=N-1-j, \  j \neq (N-1)/2, \\
(-1)^{N-1-j-k} \mathrm{sgn}(N-1-2j), & \mathrm{for} \ \mathrm{min}(j,N-1-j) < k < \mathrm{max}(j,N-1-j), \\
0, & \mathrm{otherwise}.
\end{array}
\right.
\end{eqnarray}
\begin{eqnarray}
\[ \overline{K}^{\mathrm{tr}}(r,s) \]_{j}^{k}
=s^{j-k} \times
 \left\{
\begin{array}{cc}
-r^{-1}, & \mathrm{for} \ j \le k, \ j+k=N-1, \\
-r, & \mathrm{for} \ j>k, \ j+k=N-1, \\
(-1)^{N-1-j-k}\mathrm{sgn}(N-1-2j)(r-r^{-1}), &
\mathrm{for} \ j \le k<N-1-j, \\
(-1)^{N-1-j-k}\mathrm{sgn}(N-1-2j)(r-r^{-1}),
&
\mathrm{for} \ N-1-j<k<j, \\
0 & \mathrm{otherwise}.
\end{array}
\right.
\end{eqnarray}
\hspace*{\fill}    \sq
\end{proposition}

\subsection{$N=3$}
Let us compare the
$K$-matrix obtained as a result of the degeneration and 
deBaxterization from the elliptic $K$-operator
with the full solution of the reflection equation
in the case $N=3$.
One finds that the full constant $K$-matrix is given by
\begin{eqnarray}
 K^{\mathrm{tr}} 
=\left(
\begin{array}{ccc}
d_1+d_5 & 0 & d_3 \\
d_4 & d_5 & d_6 \\
-d_7 & 0 & 0
\end{array}
\right),
\end{eqnarray}
where the solution manifold $\mathcal{S}$ of the
parameters $d_j, \ j=1,3,4,5,6,7$
is given by the Segre threefold
\begin{align}
\mathcal{S}=\{(d_1, d_3, d_4, d_5, d_6, d_7) \in \mathbf{P}^5(\mathbf{C}) \ | \
d_1 d_5-d_3 d_7=0, \ d_1 d_6-d_3 d_4=0, \ d_4 d_5-d_6 d_7=0 \}.
\end{align}
The points in the Segre threefold can be paramterized
by $\mathbf{P}^1(\mathbf{C}) \times \mathbf{P}^2(\mathbf{C})$ via
the map $\psi$
\begin{align}
&\psi:  \mathcal{U}=\mathbf{P}^1(\mathbf{C}) \times \mathbf{P}^2(\mathbf{C})
\longrightarrow \mathcal{S}, \\
&\psi((D_1, D_2) \times (E_1, E_2, E_3))
=(D_1 E_1, D_2 E_1, D_1 E_3, D_2 E_2, D_2 E_3, D_1 E_2).
\end{align}
The full constant $K$-matrix satisfies the generalized Hecke relation
\begin{equation}
(K^{\mathrm{tr}})^2-(d_1+d_5) K^{\mathrm{tr}}+d_3 d_7=0.
\end{equation}
The $K$-matrix obtained from the elliptic $K$-operator for $N=3$
\begin{eqnarray}
 K^{\mathrm{tr}}(r,s)
=\left(
\begin{array}{ccc}
r-r^{-1} &  0 & -s^2 r \\
s^{-1}(r-r^{-1}) & -r^{-1} & s(r^{-1}-r) \\
-s^{-2} r^{-1} & 0 & 0
\end{array}
\right),
\end{eqnarray}
lives on a submanifold $\mathcal{V}$ of the projective space $\mathcal{U}$
parametrizing the Segre threefold $\mathcal{S}$
\begin{align}
&\mathcal{V}=\{(-s^{-1},s) \times (-rs, -r^{-1} s^{-1}, r^{-1}-r) \in \mathcal{U} \} \nonumber \\
&=\{(D_1, D_2) \times (E_1, E_2, E_3) \in \mathbf{P}^1(\mathbf{C}) \times \mathbf{P}^2(\mathbf{C}) \ | \ D_1 D_2+1=0, \ E_1 E_2-1=0, \ D_1 E_1+D_2 E_2+E_3=0 \}.
\end{align}

We make some comments.
The full rational constant $K$-matrix can be Baxterized to give
the spectral parameter dependent $K$-matrix
\begin{align}
 K^{\mathrm{tr}}(z)&=dz^2+(d_1+d_5)z^4+(1-z^4)K^{\mathrm{tr}}
\nonumber \\
&=\left(
\begin{array}{ccc}
d_1+d_5+dz^2 &  0 & d_3(1-z^4) \\
d_4(1-z^4) & d_5+dz^2+d_1 z^4 & d_6(1-z^4) \\
-d_7(1-z^4) & 0 & dz^2+(d_1+d_5) z^4
\end{array}
\right).
\end{align}
Including the parameter $d$ which appears in the Baxterization,
the solution manifold
$\mathbf{P}^1(\mathbf{C})
\times \mathbf{P}^2(\mathbf{C})$
is lifted up to 
$\mathbf{C} \times
\mathbf{P}^1(\mathbf{C})
\times \mathbf{P}^2(\mathbf{C})$.
However, this $K$-matrix is not the full solution to the
reflection equation of the CG-$R$-matrix
with spectral parameter.
There is another solution whose solution manifold is parametrized by
$\mathbf{P}^1(\mathbf{C}) \times \mathbf{P}^1(\mathbf{C})
\times \mathbf{P}^2(\mathbf{C})$ which does not seem to be obtained at least
from a simple Baxterization scheme.
This phenomena is not observed in the standard $R$-matrix of
$U_q(\widehat{sl_3})$  \cite{Ya}.
See \cite{MY} for more details about the full solution space of $N=3$
CG-$R$-matrix.

\section{Rational representation}
For the special case $t=t_0=t_n=1$ of the type $C$ affine Hecke algebra,
one can construct another nonstandard representation
from quantum integrable models.
We call this as nonstandard rational representation
since this comes from a class of quantum integrable models of rational type.
As again, we first review how the rational Jordanian $R$-matrix \cite{EH}
is obtained from the $R$-operator.
Then we construct its associated $K$-matrix.
\subsection{Jordanian $R$-matrix}
The Jordanian $R$-matrix
can be obtained from the SU-$R$-operator as follows.
First, we degenerate the trigonometric $R$-operator furthermore.
Namely, we take the rational limit,
replacing functions $\mathrm{sin}(\pi \mathrm{i} \zeta)$ of $\zeta$ to $\zeta$.
\begin{equation}
\widetilde{\mathcal{R}}_{12}^{\mathrm{ra}}(\lambda)
:=\widetilde{\mathcal{R}}^{\mathrm{tr}}_{12}
(\lambda)|_{\mathrm{sin}(\pi \mathrm{i} \zeta) \rightarrow \zeta}.
\end{equation}
The action of the rational $R$-operator on the function 
$f(\zeta_1,\zeta_2)$ is given by
\begin{align}
\widetilde{\mathcal{R}}_{12}^{\mathrm{ra}}(\lambda)f(\zeta_1,\zeta_2)
=&\frac{\zeta_1-\zeta_2-2 \beta+\lambda}{(\zeta_1-\zeta_2-2 \beta)\lambda}f(\zeta_2,\zeta_1)
-\frac{\zeta_1-\zeta_2-2 \beta+\kappa}{(\zeta_1-\zeta_2-2 \beta) \kappa}
f(\zeta_1-2 \beta,\zeta_2+2 \beta).
\end{align}
Next, one takes the finite-dimensional representation.
We restrict the space of functions $\mathcal{M}$ to the finite-dimensional
subspace $V_N^{\mathrm{ra}} \otimes V_N^{\mathrm{ra}}$
where
\begin{equation}
V_N^{\mathrm{ra}}=\mathrm{Span}
\Big\{
\psi_k(\zeta)=\zeta^k \ \Big| \
k=0,1,\cdots,N-1 
\Big\}.
\end{equation}
Calculating the matrix elements
of the rational $R$-operator explicitly
using $\psi_k(\zeta), \ k=0,1,\cdots,N-1$ as the basis,
one gets \cite{EH,EH2}
\begin{align}
\widetilde{\mathcal{R}}_{12}^{\mathrm{ra}}(\lambda) 
\psi_i(\zeta_1) \psi_j(\zeta_2)
=\sum_{k,l=0}^{N-1}
\[ R^{\mathrm{ra}}(\lambda,\kappa,h) \]_{ij}^{kl} \psi_k(\zeta_1) \psi_l(\zeta_2),
\end{align}
with \\
\begin{align}
&\[ R^{\mathrm{ra}}(\lambda,\kappa,h) \]_{ij}^{kl}
=\frac{1}{\lambda} \delta_{il}\delta_{jk} \nonumber \\
+&
(-1)^{j-l} h^{i+j-k-l}
\Bigg\{
\left(
\begin{array}{c}
i \\
k
\end{array}
\right)
\left(
\begin{array}{c}
j \\
l
\end{array}
\right)
-\frac{\kappa}{h}
\sum_m
(-1)^{m-k}
\left(
\begin{array}{c}
i \\
m
\end{array}
\right)
\left(
\begin{array}{c}
j+m-k-1 \\
l
\end{array}
\right)
\epsilon(j,m,k)
\Bigg\},
\label{spectral-jordanian}
\end{align}
where $h=-2 \beta$.
The $R$-matrix $R^{\mathrm{ra}}(\lambda,\kappa,h)$
is called the Jordanian $R$-matrix.
From its construction from the SU-$R$-operator,
the Jordanian $R$-matrix satisfies the Yang-Baxter relation
\begin{equation}
R_{12}(\lambda)R_{13}(\lambda+\mu)R_{23}(\mu)=
R_{23}(\mu)R_{13}(\lambda+\mu)R_{12}(\lambda). \label{yang-baxter-rational}
\end{equation}
We now apply the following lemma to extract the constant $R$-matrix.
\begin{lemma} \label{gettingridoflemma-rational}
Let an $R$-matrix $R(\lambda)$ satisfying the Yang-Baxter relation
\eqref{yang-baxter-rational}
can be decomposed as
\begin{equation}
R(\lambda)=f(\lambda)(P+g \lambda R), \ f(\lambda) \not\equiv 0,
\label{decomposition-ver2-rational}
\end{equation}
and the $\lambda$-independent $R$ satisfies the degenerate Hecke relation
\begin{equation}
\check{R}^2-I=0, \label{hecke-ver2-rational}
\end{equation}
where $\check{R}=RP$.
The $\lambda$-independent $R$ satisfies the Yang-Baxter relation
without spectral parameter
\begin{equation}
R_{12}R_{13}R_{23}=
R_{23}R_{13}R_{12}. \label{ybr-ver2-rational}
\end{equation}
\hspace*{\fill}    \sq
\end{lemma}
We now apply this lemma.
One easily sees that the $R$-matrix $R^{\mathrm{ra}}(\lambda,\kappa,h)$
\eqref{spectral-jordanian} can be decomposed as
\begin{equation}
R^{\mathrm{ra}}(\lambda,\kappa,h)
=\frac{1}{\lambda}(P+\lambda R^{\mathrm{ra}}(\kappa,h)),
\label{explicitdecomposition-rational}
\end{equation}
with the $\lambda$-independent $R$-matrix given by
\begin{align}
&\[ R^{\mathrm{ra}}(\kappa,h) \]_{ij}^{kl} \nonumber \\
=&(-1)^{j-l} h^{i+j-k-l}
\Bigg\{
\left(
\begin{array}{c}
i \\
k
\end{array}
\right)
\left(
\begin{array}{c}
j \\
l
\end{array}
\right)
-\frac{\kappa}{h}
\sum_m
(-1)^{m-k}
\left(
\begin{array}{c}
i \\
m
\end{array}
\right)
\left(
\begin{array}{c}
j+m-k-1 \\
l
\end{array}
\right)
\epsilon(j,m,k)
\Bigg\}.
\end{align}
Since this $R$-matrix can be shown to satisfy the degenerate Hecke relation
\eqref{hecke-ver2-rational}, we can apply the lemma and find that the 
$R$-matrix $R^{\mathrm{ra}}(\kappa,h)$ satisfies the Yang-Baxter relation
\eqref{ybr-ver2-rational} which is equivalent to the braid relation
\begin{equation}
\check{R}_{12}\check{R}_{23}\check{R}_{12}
=\check{R}_{23}\check{R}_{12}\check{R}_{23}. \label{braid-ver2-rational}
\end{equation}
This braid relation \eqref{braid-ver2-rational}
together with the Hecke relation
\eqref{hecke-ver2-rational} can be identified with the defining relations
\eqref{relfour} and \eqref{reltwo} of type $C$ affine Hecke algebra,
hence we can use the Jordanian $R$-matrix
multiplied by the permutation matrix
$\check{R}^{\mathrm{ra}}(\kappa,h)$
as a representation for the generators
$T_1,\cdots,T_{n-1}$ of the special case $t=t_0=t_n=1$ of
the affine Hecke algebra of type $C$:
\begin{equation}
T_j=\check{R}_j^{\mathrm{ra}}(\kappa,h), \ j=1,\cdots,n-1.
\label{forthone1-rational}
\end{equation}

\subsection{Jordanian $K$-matrix}
Now we calculate the $K$-matrix of the reflection equation
of the Jordanian $R$-matrix from the HK-$K$-operator.
First, one takes furthermore the
rational limit of the trigonometric $K$-operator
\begin{equation}
\widetilde{\mathcal{K}}^{\mathrm{ra}}(\lambda)
:=\widetilde{\mathcal{K}}^{\mathrm{tr}}(\lambda)
|_{\mathrm{sin}(\pi \mathrm{i} \zeta) \to \zeta}.
\end{equation}
The action of the rational $K$-operator
on the function 
$f(\zeta)$ is given by
\begin{align}
\widetilde{\mathcal{K}}^{\mathrm{ra}}(\lambda)f(\zeta)
=&\frac{\zeta+g/2+\nu}{(\zeta+g/2)\nu}f(-\zeta-g)
-\frac{\zeta+g/2+2 \lambda}{2 (\zeta+g/2)\lambda}
f(\zeta),
\end{align}
where $g=-2 \gamma$.
Next, one takes the finite-dimensional representation.
We restrict the space of functions $\mathcal{M}$ to the finite-dimensional
subspace $V_N^{\mathrm{ra}}$.
The matrix elements of the rational $K$-operator can be calculated explicitly
using $\psi_k(\zeta), \ k=0,1,\cdots,N-1$ as the basis. One gets the following:
\begin{proposition}
The matrix elements
$\[ K^{\mathrm{ra}}(\lambda,\nu,g) \]_{j}^{k}$
of the rational $K$-operator
$\widetilde{\mathcal{K}}^{\mathrm{ra}}(\lambda)$
in the basis $\psi_k(\zeta), \ k=0,1,\cdots,N-1$
\begin{align}
\widetilde{\mathcal{K}}^{\mathrm{ra}}(\lambda) 
\psi_j(\zeta)
=\sum_{k=0}^{N-1}
\[ K^{\mathrm{ra}}(\lambda,\nu,g) \]_{j}^{k} \psi_k(\zeta),
\end{align}
is given by 
\begin{eqnarray}
\[ K^{\mathrm{ra}}(\lambda,\nu,g) \]_{j}^{k}
=-\frac{1}{2 \lambda} \delta_{jk}+\frac{(-1)^j}{\nu}
\left(
\begin{array}{c}
j \\
k
\end{array}
\right)
g^{j-k}
+2 \sum_{0 \le l<j} (-1)^{j-l}
\left(
\begin{array}{c}
j-l-1 \\
k-l
\end{array}
\right)
g^{j-k-1}.
\label{spectral-jordanian-k}
\end{eqnarray}
\hspace*{\fill}    \sq
\end{proposition}
{\it Proof}. \\
We act $\widetilde{\mathcal{K}}^{\mathrm{ra}}(\lambda)$ on $\psi_j(\zeta)$
\begin{align}
\widetilde{\mathcal{K}}^{\mathrm{ra}}(\lambda) \psi_j(\zeta)
=&-\frac{\zeta^j}{2 \lambda}+\frac{1}{\nu}(-\zeta-k)^j
+\frac{ \{-(\zeta+k/2)-k/2 \}^j-\{(\zeta+k/2)-k/2 \}^j}{\zeta+k/2}
\nonumber \\
=&-\frac{\zeta^j}{2 \lambda}+\frac{1}{\nu} \sum_k (-1)^j
\left(
\begin{array}{c}
j \\
k
\end{array}
\right)
g^{j-k} \zeta^k
+2 \sum_{k, 0 \le l<j}
(-1)^{j-l}
\left(
\begin{array}{c}
j-l-1 \\
k-l
\end{array}
\right)
g^{j-k-1}
\zeta^k.
\end{align}
\hspace*{\fill}    \wsq \\
The following proposition follows from the construction procedure from the $K$-operator.
\begin{proposition}
The matrix $K^{\mathrm{ra}}(\lambda,\nu,g)$
is a solution to the reflection equation
of the Jordanian rational $R$-matrix $R^{\mathrm{ra}}(\lambda,\kappa,h)$.
\hspace*{\fill}    \sq
\end{proposition}
{\it Proof}. \\
This follows from the fact that the $K$-matrix
is constructed as a degeneration from the HK-$K$-operator,
which is a solution to the reflection equation of the
SU-$R$-operator. \hspace*{\fill}    \wsq \\
So far, we found a solution $K^{\mathrm{ra}}(\lambda,\nu,g)$
to the reflection equation of the Jordanian $R$-matrix
$R^{\mathrm{ra}}(\lambda,\kappa,h)$ with spectral parameter
\begin{equation}
R_{12}(\lambda-\mu)
K_{1}(\lambda)
R_{21}(\lambda+\mu)
K_{2}(\mu)
=
K_{2}(\mu)
R_{12}(\lambda+\mu)
K_{1}(\lambda)
R_{21}(\lambda-\mu). \label{reflectionforlemma-rational}
\end{equation}
To extract representations for the generators $T_0$ and $T_n$,
we use the following lemma.
\begin{lemma} \label{auxlemmarational}
Let $R(\lambda)$ be an $R$-matrix which satisfies the properties in
Lemma \ref{gettingridoflemma-rational}.
If the corresponding $K$-matrix $K(\lambda)$ satisfying
the reflection equation
\eqref{reflectionforlemma-rational}
can be decomposed as
\begin{equation}
K(\lambda)=a(\lambda)(I+b\lambda K), \ a(\lambda) \not\equiv 0,
\label{decompositionrational}
\end{equation}
and the $\lambda$-independent $K$ satisfies the Hecke relation
\begin{equation}
K^2-I=0, \label{heckerational}
\end{equation}
then the $\lambda$-independent $K$ satisfies the reflection equations
without spectral parameter
\begin{align}
&\check{R}_{12}K_{1}\check{R}_{12}K_{1}
=K_{1}\check{R}_{12}K_{1}\check{R}_{12}, \label{whatweneed1rational} \\
&\check{R}_{12}K_{2}\check{R}_{12}K_{2}
=K_{2}\check{R}_{12}K_{2}\check{R}_{12}. \label{whatweneed2rational}
\end{align}
\hspace*{\fill}    \sq
\end{lemma}
{\it Proof}. \\
Multiplying both sides of the spectral parameter dependent reflection equation
\eqref{reflectionforlemma-rational} by the permutation operators $P$
and inserting the decomposition relations \eqref{decomposition-ver2-rational},
\eqref{decompositionrational}
and Hecke relations \eqref{hecke-ver2-rational},
\eqref{heckerational}
into it,
one finds the coefficients of the terms $K \check{R}$
and $\check{R} K$ have the following form
\begin{align}
bg( \mu + \lambda )(\lambda-\mu)+bg(\mu-\lambda)(\lambda+\mu),
\end{align}
which is obviously zero.
The remaining equations are nothing but
the reflection equations without spectral parameters
\eqref{whatweneed1rational}, \eqref{whatweneed2rational}.
\hspace*{\fill}    \wsq \\
We now apply this lemma to extract the constant $K$-matrix.

\begin{proposition} \label{auxproprational}
The $K$-matrix $K^{\mathrm{ra}}(\nu,g)$
whose matrix elements are explicitly given by
\begin{eqnarray}
\[ K^{\mathrm{ra}}(\nu,g) \]_{j}^{k}
=(-1)^j
\left(
\begin{array}{c}
j \\
k
\end{array}
\right)
g^{j-k}
+2 \nu \sum_{0 \le l<j} (-1)^{j-l}
\left(
\begin{array}{c}
j-l-1 \\
k-l
\end{array}
\right)
g^{j-k-1},
\label{forproprationalrational}
\end{eqnarray}
is a solution to the constant reflection equation
of the Jordanian $R$-matrix $R^{\mathrm{ra}}(\nu,g)$.
\hspace*{\fill}    \sq
\end{proposition}
{\it Proof}. \\
The $K$-matrix $K^{\mathrm{ra}}(\lambda,\nu,g)$ satisfying the
reflection equation \eqref{reflectionforlemma-rational}
can be decomposed as
\begin{align}
K^{\mathrm{ra}}(\lambda,\nu,g)=-\frac{1}{2 \lambda}
\Bigg(
I-\frac{2 \lambda}{\nu} K^{\mathrm{ra}}(\nu,g)
\Bigg),
\label{propdecompositionrational}
\end{align}
with the $\lambda$-dependent $K$-matrix
$K^{\mathrm{ra}}(\nu,g)$ explicitly given by \eqref{forproprationalrational}.
The Hecke relation \eqref{heckerational}
which one needs to apply Lemma \ref{auxlemmarational}
follows by comparing the expression
\begin{align}
K^{\mathrm{ra}}(\lambda,\nu,g)K^{\mathrm{ra}}(-\lambda,\nu,g)
=\Bigg( \frac{1}{4 \lambda^2}+\frac{1}{\nu^2} \Bigg)I,
\end{align}
which can be calculated using the operator expression for
$K^{\mathrm{ra}}(\lambda,\nu,g)$, and comparing with another expression
\begin{align}
K^{\mathrm{ra}}(\lambda,\nu,g)K^{\mathrm{ra}}(-\lambda,\nu,g)
=\frac{1}{4 \lambda^2}I+\frac{1}{\nu^2}K^2,
\end{align}
obtained from the decomposition
\eqref{propdecompositionrational}. \\
We have shown all the conditions the $K$-matrices should satisfy
to apply Lemma \ref{auxlemmarational},
and the proposition follows.
\hspace*{\fill}    \wsq \\
Lemma \ref{auxlemmarational} and Proposition \ref{auxproprational}
shows that the $K$-matrix
$K^{\mathrm{ra}}(\nu,g)$ satisfies
the reflection equations
\eqref{whatweneed1rational}, \eqref{whatweneed2rational}
and the Hecke relation \eqref{heckerational}
which can be identified with the
defining relations \eqref{relfour}, \eqref{relsix},
\eqref{relone} and \eqref{relthree}
of type $C$ affine Hecke algebra,
hence one obtains a representation for $T_0$ and $T_n$
in terms of the $K$-matrix:
\begin{align}
T_0=K_1^{\mathrm{ra}}(\nu_0,g_0), \label{forthone2-rational} \\
T_n=K_n^{\mathrm{ra}}(\nu_n,g_n). \label{forthone3-rational}
\end{align}
The representation for the generators in terms
of the constant Jordanian $R$-matrix and
its associated $K$-matrix
\eqref{forthone1-rational},
\eqref{forthone2-rational} and \eqref{forthone3-rational}
proves Theorem \ref{thtwo}.
Namely, one has a representation for the affine Hecke algebra $H_n(1,1,1)$.

\subsection{$N=3$}
For $N=3$, we find that the full $K$-matrix of the constant reflection equation
for the Jordanian $K$-matrix is given as
\begin{eqnarray}
 K^{\mathrm{ra}} 
=\left(
\begin{array}{ccc}
c_1 & c_2 & c_3 \\
0 & c_5 & c_6 \\
0 & 0 & c_9
\end{array}
\right).
\end{eqnarray}
Here the parameters $c_1, c_2, c_3, c_5, c_6, c_9$
live on the following solution manifold
\begin{align}
&c_2 c_6+c_3(c_1-c_5)=0, \\
&c_2(c_1-c_9)=0, \\
&(c_1-c_5)(c_1-c_9)=0.
\end{align}
One can show that the full constant $K$-matrix
satisfies the generalized Hecke relation
\begin{equation}
(K^{\mathrm{ra}})^2-(c_5+c_9)K^{\mathrm{ra}}+c_5 c_9=0. \label{fullhecke}
\end{equation}
Analyzing the solution manifold,
one finds the constant $K$-matrix
can be furthermore divided into two types
\begin{eqnarray}
 K_{\mathrm{I}}^{\mathrm{ra}} 
=\left(
\begin{array}{ccc}
c_1 & c_2 & \alpha c_2 \\
0 & c_5 & \alpha(c_5-c_1) \\
0 & 0 & c_1
\end{array}
\right), \ \alpha \in \mathbf{C},
\end{eqnarray}
and
\begin{eqnarray}
 K_{\mathrm{II}}^{\mathrm{ra}} 
=\left(
\begin{array}{ccc}
c_1 & 0 & c_3 \\
0 & c_1 & c_6 \\
0 & 0 & c_9
\end{array}
\right).
\end{eqnarray}
The solution manifold $\mathcal{B}$ of the
first solution $K_{\mathrm{I}}^{\mathrm{ra}}$ is the projective space
$\mathbf{P}^3(\mathbf{C})$
\begin{align}
\mathcal{A}=\{(c_1, c_2, c_5, \alpha) \in \mathbf{P}^3(\mathbf{C})\}.
\end{align}
The solution obtained as a degeneration from the elliptic $K$-operator
\begin{eqnarray}
 K^{\mathrm{ra}} 
=\left(
\begin{array}{ccc}
1 & -g-2 \nu & g^2+2g \nu \\
0 & -1 & 2g \\
0 & 0 & 1
\end{array}
\right),
\end{eqnarray}
multiplied by an overall factor
lives on a hyperplane $\mathcal{B}$ of the solution manifold $\mathcal{A}$
of the first solution $K_{\mathrm{I}}^{\mathrm{ra}}$
\begin{align}
\mathcal{B}=\{(c_1, c_2, c_5, \alpha) \in \mathbf{P}^3(\mathbf{C}) \ | \
c_5=-c_1 \}.
\end{align}
The solution obatined as a degeneration and deBaxterization
from the $K$-operator
can construct representations only for the special case $H_n(1,1,1)$.
On the other hand, the Hecke relation \eqref{fullhecke} shows that
the full solution can construct representations for $H_n(1,t_0,t_n)$.
We finally remark that the full constant $K$-matrix can be Baxterized
to give the spectral parameter dependent $K$-matrix
\begin{align}
 K^{\mathrm{ra}}(\lambda)&=c-\lambda(c_5+c_9)
+2 \lambda K^{\mathrm{ra}}
\nonumber \\
&=\left(
\begin{array}{ccc}
c+\lambda(2c_1-c_5-c_9) &  2 \lambda c_2 & 2 \lambda c_3 \\
0 & c+\lambda(c_5-c_9) & 2 \lambda c_6 \\
0 & 0 & c+\lambda(c_9-c_5)
\end{array}
\right).
\end{align}
\section{Discussion}
In this paper, we constructed explicit nonstandard representations of
type $C$ affine Hecke algebra.
Concretely realizing representations of affine Hecke algebra
is not an easy problem. We can approach this problem by
using the power of quantum integrable models.
For type $C$ affine Hecke algebra,
we achieved this by using two classes of quantum integrable models
under the reflecting boundary condition.
The nonstandard Cremmer-Gervais $R$-matrix serves as representations
for the generators of type $A$ affine Hecke algebra since it satisfies
the Hecke relation as well as the Yang-Baxter relation.
To construct nonstandard representations for type $C$ is equivalent
to finding solutions of reflection equation under the Hecke relation
(see \cite{No,Sa} for standard representations of
type $C$ affine Hecke algebra based on standard $R$ and $K$-matrices).

We constructed them by taking appropriate degeneration and deBaxterization
of the Hikami-Komori elliptic $K$-operator.
We also constructed another representation for a special case of type
$C$ affine Hecke algebra from the rational Jordanian $R$-matrix and
its associated $K$-matrix, also achieved by the degeneration and
deBaxterization scheme from the $R$-operator and the $K$-operator.
The degeneration procedure in \cite{EH} and in this paper seems to show
a systematic way of constructing representations of affine Hecke algebras
not found yet.
Finding nonstandard representations do not have to seem recipes.
However, starting from the infinite-dimensional operators
and taking finite-dimensional representations can yield
nonstandard representations,
and the advantage of starting from infinite-dimensional operators
is that proving relations at the level of operators
are much easier than at the level of finite-dimensional representations.
It may be worth investigating affine Hecke algebras
associated with other root systems in this way for example,
and is also worthwhile to investigate affine Hecke algebras of type $C$
from the point of view of
boundary quantum group \cite{Ne}.
It may also be interesting to use the results in this paper
to formulate and study boundary (type $C$)
analogue of the nilpotency indices of the $R$-matrices (type $A$) \cite{GG,Jo},
or to relate other integrable systems such as the classical top
(see \cite{LOZ} for an example of relating nonstandard $R$-matrix with
integrable tops).
\section*{Acknowledgments}
This work was partially supported by
grants-in-aid for Scientific Research (C) No. 24540393 and
for Young Scientists (B) No. 25800223.
%
%
%
%

\end{document}